\documentclass[11pt]{article}
\usepackage{fullpage,amssymb}
\title{{\bf Vertex tensor category structure on a category of
Kazhdan--Lusztig}}
\author{Lin Zhang}
\date{}

\newtheorem{rema}{Remark}[section]
\newtheorem{propo}[rema]{Proposition}
\newtheorem{theo}[rema]{Theorem}
\newtheorem{defi}[rema]{Definition}
\newtheorem{lemma}[rema]{Lemma}

\newcommand{\nno}{\nonumber}
\newcommand{\lbar}{\bigg\vert}
\newcommand{\pf}{{\it Proof}\hspace{2ex}}
\newcommand{\epf}{\hspace{2em}$\square$}
\newcommand{\epfv}{\hspace{1em}$\square$\vspace{1em}}

\newcommand{\res}{\mbox{\rm Res}}
\newcommand{\wt}{\mbox{\rm wt}\ }
\newcommand{\ob}{{\rm ob}\,}
\renewcommand{\hom}{\mbox{\rm Hom}}

 \makeatletter
\newlength{\@pxlwd} \newlength{\@rulewd} \newlength{\@pxlht}
\catcode`.=\active \catcode`B=\active \catcode`:=\active \catcode`|=\active
\def\sprite#1(#2,#3)[#4,#5]{
   \edef\@sprbox{\expandafter\@cdr\string#1\@nil @box}
   \expandafter\newsavebox\csname\@sprbox\endcsname
   \edef#1{\expandafter\usebox\csname\@sprbox\endcsname}
   \expandafter\setbox\csname\@sprbox\endcsname =\hbox\bgroup
   \vbox\bgroup
  \catcode`.=\active\catcode`B=\active\catcode`:=\active\catcode`|=\active
      \@pxlwd=#4 \divide\@pxlwd by #3 \@rulewd=\@pxlwd
      \@pxlht=#5 \divide\@pxlht by #2
      \def .{\hskip \@pxlwd \ignorespaces}
      \def B{\@ifnextchar B{\advance\@rulewd by \@pxlwd}{\vrule
         height \@pxlht width \@rulewd depth 0 pt \@rulewd=\@pxlwd}}
      \def :{\hbox\bgroup\vrule height \@pxlht width 0pt depth
0pt\ignorespaces}
      \def |{\vrule height \@pxlht width 0pt depth 0pt\egroup
         \prevdepth= -1000 pt}
   }
\def\endsprite{\egroup\egroup}
\catcode`.=12 \catcode`B=11 \catcode`:=12 \catcode`|=12\relax
\makeatother

\def\hboxtr{\FormOfHboxtr} 
\sprite{\FormOfHboxtr}(25,25)[0.5 em, 1.2 ex] 

:BBBBBBBBBBBBBBBBBBBBBBBBB |
:BB......................B |
:B.B.....................B |
:B..B....................B |
:B...B...................B |
:B....B..................B |
:B.....B.................B |
:B......B................B |
:B.......B...............B |
:B........B..............B |
:B.........B.............B |
:B..........B............B |
:B...........B...........B |
:B............B..........B |
:B.............B.........B |
:B..............B........B |
:B...............B.......B |
:B................B......B |
:B.................B.....B |
:B..................B....B |
:B...................B...B |
:B....................B..B |
:B.....................B.B |
:B......................BB |
:BBBBBBBBBBBBBBBBBBBBBBBBB |

\endsprite

\sprite{\FormOfShboxtr}(25,25)[0.3 em, 0.72 ex] 

:BBBBBBBBBBBBBBBBBBBBBBBBB |
:BB......................B |
:B.B.....................B |
:B..B....................B |
:B...B...................B |
:B....B..................B |
:B.....B.................B |
:B......B................B |
:B.......B...............B |
:B........B..............B |
:B.........B.............B |
:B..........B............B |
:B...........B...........B |
:B............B..........B |
:B.............B.........B |
:B..............B........B |
:B...............B.......B |
:B................B......B |
:B.................B.....B |
:B..................B....B |
:B...................B...B |
:B....................B..B |
:B.....................B.B |
:B......................BB |
:BBBBBBBBBBBBBBBBBBBBBBBBB |

\endsprite

\newcommand{\fg}{{\mathfrak g}}
\newcommand{\fh}{{\mathfrak h}}
\newcommand{\fn}{{\mathfrak n}}
\newcommand{\C}{{\mathbb C}}
\newcommand{\N}{{\mathbb N}}
\newcommand{\Q}{{\mathbb Q}}
\newcommand{\Z}{{\mathbb Z}}
\newcommand{\bk}{{\bf k}}

\begin{document}
\maketitle

\begin{abstract}
We incorporate a category of certain modules for an affine Lie
algebra, of a certain fixed non-positive-integral level, considered by Kazhdan and
Lusztig, into the representation theory of vertex operator algebras, by using
the logarithmic tensor product theory for generalized modules for a vertex
operator algebra developed by Huang,
Lepowsky and the author. We do this by proving that the conditions for
applying this general logarithmic tensor product theory hold. 
As a consequence, we prove that this category
has a natural vertex tensor category structure, and in particular
we obtain a new, vertex-algebraic, construction of the natural 
associativity isomorphisms and proof of their properties.
\end{abstract}

\renewcommand{\theequation}{\thesection.\arabic{equation}}
\renewcommand{\therema}{\thesection.\arabic{rema}}

\setcounter{equation}{0} \setcounter{rema}{0}
\section{Introduction}
In a series of papers \cite{KL2}--\cite{KL5} (see also \cite{KL1}),
Kazhdan and Lusztig constructed a braided tensor category structure on
a category ${\cal O}_\kappa$ of certain modules for an affine Lie
algebra of a fixed level $\kappa-h$, $h$ the dual Coxeter number of
the Lie algebra and $\kappa$ a complex number not in $\Q_{\geq 0}$,
the set of nonnegative rational numbers,
and showed that this braided tensor category is equivalent to a tensor
category of modules for a quantum group constructed from the same Lie
algebra. 
The most interesting cases are those of the allowable negative levels.
Their construction of the tensor product uses ideas in the
important work \cite{MS} by Moore and Seiberg. While the category in
\cite{MS} consists essentially of the integrable highest weight
modules for an affine Lie algebra of a
fixed positive integral level rather than modules of a fixed negative
level, for example, it
was first discovered by Moore and Seiberg in that paper that this positive-level
category has a braided tensor category structure; however, this was based on the
strong assumption of the axioms for conformal field theory and in particular, the
assumption of the existence of the ``operator product expansion'' for 
``intertwining operators.''

Suitable modules for affine Lie algebras give rise to an important
family of vertex operator algebras and their modules. In
\cite{tensorK}--\cite{tensor3} and \cite{tensor4}, Huang and Lepowsky
developed a substantial tensor product theory for modules for a
``rational'' vertex operator algebra, under certain conditions.  As
one of the applications of this theory, they proved in \cite{HLaffine}
that the conditions required for this theory are satisfied for the
module category of the vertex operator algebra constructed
{}from the category of integrable highest weight modules for an affine
Lie algebra of a fixed positive integral level. As a result they
directly constructed a braided tensor category structure, and further,
a ``vertex tensor category'' structure, on this category.

It has been expected that the category ${\cal O}_\kappa$ considered by
Kazhdan and Lusztig should also be covered by a suitable generalization
of the tensor product theory developed by Huang and Lepowsky, even
though $\kappa-h$ cannot be positive integral.
However, this category is very different from the one associated with
the positive integral
level case. For example, the objects of this category in
general are only direct sums
of generalized eigenspaces, rather than eigenspaces, for the operator
$L_0$ that defines the conformal weights. Therefore, such a
generalization, if it exists, should include categories of these more
general modules. Recently, such a generalization has been achieved in
\cite{HLZ1} and \cite{HLZ2} by Huang, Lepowsky and the author. A fundamental
subtlety in this generalization is that the corresponding intertwining
operators involve logarithms of the variables.  The questions are
now whether the category of Kazhdan and Lusztig satisfies the required
conditions for the generalized tensor product theory of \cite{HLZ1}
and \cite{HLZ2}, and if the
answer is affirmative, whether the resulting tensor product
construction is equivalent to the original construction given by
Kazhdan and Lusztig.

In this paper, we prove that category ${\cal O}_\kappa$ indeed
satisfies all necessary conditions in \cite{HLZ1} and \cite{HLZ2}.  We
establish an equivalent condition for the subtle ``compatibility
condition'' in the construction of the tensor product and use it to
prove that the two constructions of the tensor product functor are
identical. We then use the methods in \cite{tensor4} and \cite{H3} and
their generalizations in \cite{HLZ1} and \cite{HLZ2} to obtain a new
construction, very different from the original one by Kazhdan and
Lusztig, of the natural associativity isomorphisms.  As a result, we
incorporate the tensor category theory of ${\cal O}_\kappa$ into the theory of vertex operator
algebras and more importantly, prove that ${\cal O}_\kappa$ has a
natural vertex tensor category structure.

The contents of this paper are as follows: In Section 2 we recall the
construction of the tensor product by Kazhdan and Lusztig and some of
their results on the category ${\cal O}_\kappa$. In Section 3 we
recall from \cite{HLZ1} and \cite{HLZ2} the construction of tensor
product of generalized modules for a vertex operator algebra. In Section 4 we
prove an equivalent condition for the ``compatibility condition.''  Then in
Section 5 we first apply the general theory to the case of ${\cal
O}_\kappa$ and prove the equivalence of the two constructions of the
tensor product functor; then we show that the objects of the category ${\cal
O}_\kappa$ are $C_{1}$-cofinite and quasi-finite dimensional.
These results imply that the
conditions for applying the results in \cite{HLZ1} to the
category ${\cal O}_\kappa$ are
satisfied, and thus ${\cal O}_\kappa$ has a vertex
tensor category structure.  This paper is heavily based on the
generalized tensor product theory developed in \cite{HLZ1} and
\cite{HLZ2}.

In this paper, $\C$, $\N$ and $\Z_{+}$ are the complex numbers, 
the nonnegative integers and the positive integers, respectively.

\paragraph{Acknowledgment} 
This paper is a revised version of part of the author's
Ph.D. thesis \cite{Z}.
The author would like to thank Professor James Lepowsky and Professor
Yi-Zhi Huang for their long time encouragement and support, both
spiritually and financially, and for numerous discussions and suggestions
on this and
related works; he would also like to thank Professor Haisheng Li for
his encouragement and discussions from the very beginning. Part of the
contents of this paper was presented in seminar talks given at Stony
Brook and at Rutgers by the author in April and November, 2003. The
author would like to thank Bin Zhang for inviting him to present
this work at Stony Brook. Partial support from NSF grant DMS-0070800
is gratefully acknowledged.

\setcounter{equation}{0} \setcounter{rema}{0}
\section{Kazhdan-Lusztig's tensor product and the category
${\cal O}_\kappa$}

In this section we recall the ``double dual'' construction of tensor
product of certain modules for an affine Lie algebra of a fixed level,
given by Kazhdan and Lusztig in \cite{KL2} (see also \cite{KL1}). We
also recall from their papers the category ${\cal O}_\kappa$ for a
complex number $\kappa\notin\Q_{\geq 0}$ and the result on closedness
of tensor product on this category (see also \cite{Y}).

Let $\fg$ be a complex semisimple finite dimensional Lie algebra
equipped with a nondegenerate invariant symmetric bilinear form
$(\cdot,\cdot)$. The {\em (untwisted) affine Lie algebra associated
with $\fg$} is the vector space
\[
\hat\fg=\fg\otimes\C[t,t^{-1}]\oplus\C\bk,
\]
equipped with the bilinear bracket operations
\begin{eqnarray}
&[a\otimes t^m,b\otimes t^n]=[a,b]\otimes t^{m+n}+
m(a,b)\delta_{m+n,0}\bk,&\label{ghatbracket}\\
&[\cdot,\bk]=0=[\bk,\cdot]\label{kcentral}
\end{eqnarray}
for $a,b\in\fg$ and $m,n\in\Z$. We will also need its algebraic 
completion $\fg\otimes \C((t))\oplus\C\bk$, which satisfies the same
bracket relations (\ref{ghatbracket}) and (\ref{kcentral}), and more
generally, for $a,b\in\fg$ and $g_1, g_2\in\C((t))$,
\begin{equation}\label{ghatbracketc}
[a\otimes g_1,b\otimes g_2]=[a,b]\otimes g_1g_2+
\{g_1,g_2\}(a,b)\bk
\end{equation}
where $\{g_1,g_2\}=\res\,g_2\frac{d}{dt}g_1$, the coefficient of
$t^{-1}$ in the formal Laurent series $g_2\frac{d}{dt}g_1$.

Equipped with the $\Z$-grading
\[
\hat\fg=\coprod_{n\in\Z}\hat\fg_{(n)},
\]
where
\[
\hat\fg_{(0)}=\fg\oplus\C\bk\;\mbox{ and }\;
\hat\fg_{(n)}=\fg\otimes t^{-n}\mbox{ for }n\neq 0,
\]
$\hat\fg$ becomes a $\Z$-graded Lie algebra.  We have the following
graded subalgebras of $\hat\fg$:
\begin{eqnarray*}
\hat\fg_{(\pm)}&=&\coprod_{n>0}\fg\otimes t^{\mp n},\\
\hat\fg_{(\leq 0)}&=&\hat\fg_{(-)}\oplus\fg\oplus\C\bk.
\end{eqnarray*}

A $\hat\fg$-module $W$ is said to be of {\em level $\ell\in\C$} if
$\bk$ acts on $W$ as scalar $\ell$.
A module
$W$ for $\hat\fg$ or $\hat\fg_{(\leq 0)}$ is said to be {\em
restricted} if for any $a\in\fg$ and $w\in W$, $(a\otimes t^n)w=0$ for
$n$ sufficiently large.  Note that a restricted module $W$ for
$\hat\fg$ is naturally a module for $\fg\otimes \C((t))\oplus\C \bk$
by letting $a\otimes\sum_{n\in\Z}c_n t^n$, $a\in\fg$, $c_n\in\C$, act
on $w\in W$ as $\sum_{n\in\Z}c_n(a\otimes t^n)w$.

A $\hat\fg$-module $W$ is said to be {\em smooth} if for any $w\in W$,
there is $N\in\N$ such that for any $a_1, \dots, a_N\in \fg$,
$(a_1\otimes t)\cdots(a_N\otimes t)w=0$. By (\ref{ghatbracket}) and
the fact that $\fg=[\fg,\fg]$ it is clear that a smooth
$\hat\fg$-module must be restricted.

Let $m$ be a positive integer. First recall that the direct sum
$\fg^{\oplus m}$ of $m$ copies of $\fg$ is a semisimple Lie algebra
with nondegenerate invariant symmetric bilinear form $(\cdot,\cdot)$
given by
\[
((a_1,\cdots,a_m),(b_1,\cdots,b_m))=(a_1,b_1)+\cdots+(a_m,b_m)
\]
for $a_1, \dots, a_m, b_1, \dots, b_m\in\fg$.

Given $m$ restricted
$\hat\fg$-modules of a fixed level $\ell\in\C$, the goal is to produce
a ``tensor product'' of these modules that is also a $\hat\fg$-module
of the {\em same} level. (Note that the usual tensor product for Lie
algebra modules is a module of level $m\ell$.)  As we
recall from \cite{KL1} and \cite{KL2} below, 
this can be defined in terms of a Riemann sphere with
$m+1$ distinct points and local coordinates at these points.

Let $p_0, p_1, \dots, p_m$ be distinct points, or {\em punctures} (see
\cite{H1}), on the Riemann sphere $C=\C P^1$ and let $\varphi_s: C\to
\C P^1$ be isomorphisms such that $\varphi_s(p_s)=0$ for each of $s=0,
1,\dots,m$, that is, $\varphi_s$ is the {\em local coordinate around
$p_s$} for each $s$. We will still use $C$ for the Riemann
sphere equipped with these punctures and local coordinates.

Let $R$ denote the algebra of regular functions on $C\backslash\{p_0,
p_1, \dots, p_m\}$. Define $\{f_1,f_2\}=\res_{p_0}f_2df_1$, i.e., the
residue of the meromorphic $1$-form $f_2df_1$ on $C$ at the point
$p_0$. Then $\{\cdot,\cdot\}:R\times R\to\C$ is a bilinear form
satisfying
\[
\{f_1,f_2\}+\{f_2,f_1\}=0, \qquad
\{f_1f_2,f_3\}+\{f_2f_3,f_1\}+\{f_3f_1,f_2\}=0
\]
for all $f_1,f_2,f_3\in R$. As a result the Lie algebra $\fg\otimes R$
has a natural central extension $\Gamma_R=(\fg\otimes R)\oplus\C\bk$
with central element $\bk$ and bracket relations
\begin{equation}\label{Gammabracket}
[a\otimes f_1,b\otimes f_2]=[a,b]\otimes f_1f_2 +\{f_1,f_2\}(a,b)\bk,
\end{equation}
for $a, b\in \fg$ and $f_1, f_2\in R$.

\begin{rema}{\rm
In \cite{KL2}, $C$ is allowed to be a smooth curve with $k$ connected
components each of which is isomorphic to $\C P^1$, and the
construction would give a $\hat\fg^{\oplus k}$-module as the tensor
product. For the purpose of this paper, however, we need only 
the $k=1$ case.  }
\end{rema}

For each $s=0,1,\dots,m$, denote by $\iota_{p_s}:R\to\C((t))$ the
linear map which sends $f\in R$ to the power series expansion of
$f\circ\varphi^{-1}_s$ around $0$. Then we have Lie algebra
homomorphisms
\begin{eqnarray}\label{hom1}
&\Gamma_R\to\fg\otimes\C((t))\oplus\C\bk,&\nno\\
&a\otimes f\mapsto a\otimes\iota_{p_0}f,\quad \bk\mapsto\bk&
\end{eqnarray}
for $a\in\fg$, $f\in R$, and
\begin{eqnarray}\label{hom2}
&\Gamma_R\to(\fg\otimes\C((t)))^{\oplus m}\oplus\C\bk,&\nno\\
&a\otimes f\mapsto(a\otimes\iota_{p_1}f,\cdots,a\otimes
\iota_{p_m}f),\quad\bk\mapsto -\bk&
\end{eqnarray}
for $a\in\fg$ and $f\in R$ (cf. Remark \ref{residuethm} below).
Here we see $(\fg\otimes\C((t)))^{\oplus m}
\oplus\C\bk=\fg^{\oplus m}\otimes\C((t))\oplus\C\bk$ as the algebraic
completion of the affine Lie algebra $\widehat{\fg^{\oplus
m}}=\fg^{\oplus m}\otimes\C[t,t^{-1}]\oplus\C\bk$ 
associated with the semisimple
Lie algebra $\fg^{\oplus m}$; so in particular,
\begin{eqnarray}\label{gmbracketc}
\lefteqn{[(a_1\otimes g_1,\cdots,a_m\otimes g_m),(a'_1\otimes
g'_1,\cdots, a'_m\otimes g'_m)]}\nno\\
&&=([a_1,a'_1]\otimes g_1g'_1,\cdots,[a_m,a'_m]\otimes g_m g'_m)
+\sum_{i=1}^m(a_i,a'_i)\{g_i,g'_i\}\bk
\end{eqnarray}
for $a_1,\dots,a_m,a'_1,\dots,a'_m\in\fg$ and $g_1,\dots,g_m,g'_1,
\dots,g'_m\in\C((t))$.
\begin{rema}\label{residuethm}{\rm
Checking that (\ref{hom1}) is a Lie algebra homomorphism is
straightforward by (\ref{Gammabracket}) and
(\ref{ghatbracketc}). Formula (\ref{hom2}) gives a Lie algebra
homomorphism due to (\ref{Gammabracket}), (\ref{gmbracketc}) and the
fact that for any $f_1,f_2\in R$,
\[
\{\iota_{p_1}f_1,\iota_{p_1}f_2\}+\cdots+
\{\iota_{p_m}f_1,\iota_{p_m}f_2\}=
-\{\iota_{p_0}f_1,\iota_{p_0}f_2\},
\]
by the residue theorem. }
\end{rema}

Now let $W_1$, $W_2$, $\dots$, $W_m$ be $\hat\fg$-modules of level
$\ell$. The vector space tensor product $W_1\otimes\cdots\otimes W_m$
is naturally a module for the affine Lie algebra $\widehat{\fg^{\oplus
m}}=\fg^{\oplus m}\otimes\C[t,t^{-1}]\oplus\C\bk=(\fg\otimes
\C[t,t^{-1}])^{\oplus m}\oplus\C\bk$ by $\bk$ acting as scalar $\ell$
and
\begin{eqnarray}\label{aaawww}
\lefteqn{(a_1\otimes t^{n_1},\cdots,a_m\otimes t^{n_m})(w_{(1)}
\otimes\cdots\otimes w_{(m)})=}\nno\\
&&\hspace{1em}=(a_1\otimes t^{n_1})w_{(1)}\otimes w_{(2)}\otimes
\cdots\otimes w_{(m)}+\cdots+\nno\\
&&\hspace{2em}+w_{(1)}\otimes w_{(2)}\otimes\cdots\otimes
(a_m\otimes t^{n_m})w_{(m)}
\end{eqnarray}
for $a_1,\dots,a_m\in\fg$, $n_1,\dots,n_m\in\Z$ and $w_{(i)}\in W_i$,
$i=1,\dots,m$. This follows from the bracket relations
(\ref{gmbracketc}) where $g_i$ and $g'_i$'s are in $\C[t,t^{-1}]$.

If each $\hat\fg$-module $W_i$ is restricted, $i=1,\dots,m$, then it
is clear from
(\ref{aaawww}) that $W_1\otimes\cdots\otimes W_m$ is a restricted
$\widehat{\fg^{\oplus m}}$-module, and hence naturally a module for
$\fg^{\oplus m}\otimes\C((t))\oplus\C\bk$, satisfying
\begin{eqnarray}\label{aaawwwc}
\lefteqn{(a_1\otimes g_1,\cdots,a_m\otimes g_m)(w_{(1)}
\otimes\cdots\otimes w_{(m)})=}\nno\\
&&\hspace{1em}=(a_1\otimes g_1)w_{(1)}\otimes w_{(2)}\otimes
\cdots\otimes w_{(m)}+\cdots+\nno\\
&&\hspace{2em}+w_{(1)}\otimes w_{(2)}\otimes\cdots\otimes
(a_m\otimes g_m)w_{(m)}.
\end{eqnarray}
for $a_1,\dots,a_m\in\fg$, $g_1,\dots,g_m\in\C((t))$ and $w_{(i)}\in
W_i$, $i=1,\dots,m$. Thus by (\ref{hom2}) we have a $\Gamma_R$-module
structure on $W_1\otimes\cdots\otimes W_m$ where $\bk$ acts as the
scalar $-\ell$. The dual vector space $(W_1\otimes\cdots\otimes
W_m)^*=\hom(W_1\otimes\cdots\otimes W_m,\C)$ has an induced natural
$\Gamma_R$-module structure by
\begin{equation}\label{Liedual}
\langle\xi(\lambda),w\rangle=-\langle\lambda,\xi(w)\rangle
\end{equation}
for all $\xi\in\Gamma_R$, $\lambda\in(W_1\otimes\cdots\otimes W_m)^*$
and $w\in W_1\otimes\cdots\otimes W_m$. Here and below we use
$\langle\cdot,\cdot\rangle$ to denote the natural pairing between a
vector space and its dual.

Let $N$ be a positive integer. Let $G_N$ be the subspace of
$U(\Gamma_R)$ spanned by all products $(a_1\otimes
f_1)\cdots(a_N\otimes f_N)$ with $a_1,\dots,a_N\in\fg$ and
$f_1,\dots,f_N\in R$ satisfying $\iota_{p_0}f_i\in t\C[[t]]$ for
$i=1,\dots, N$. 
Here and below we use notation $U(L)$ for the univeral enveloping
algebra of a Lie algebra $L$.
Define $Z^N\subset(W_1\otimes\cdots\otimes W_m)^*$ to
be the annihilator of $G_N(W_1\otimes\cdots\otimes W_m)$. Then
\[
Z^N=\{\lambda\in(W_1\otimes\cdots\otimes W_m)^*\,|\, G_N\lambda=0\},
\]
and we have an increasing sequence $Z^1\subset Z^2\subset\cdots$. Let
$Z^\infty=\cup_{N\in\Z_+}Z^N$. It is clear that $Z^\infty$ is a
$\Gamma_R$-submodule of $(W_1\otimes\cdots\otimes W_m)^*$.

Define a $\fg\otimes\C((t))\oplus\C\bk$-module structure on $Z^\infty$
as follows: $\bk$ acts as scalar $\ell$; and
for $\lambda\in Z^\infty$, $a\in\fg$ and $g\in\C((t))$,
choose $N\in\N$ such that $\lambda\in Z^N$, choose $f\in R$ such that
$\iota_{p_0}f-g\in t^N\C[[t]]$, and define
\begin{equation}\label{hf}
(a\otimes g)\lambda=(a\otimes f)\lambda.
\end{equation}
It is easy to verify that this is independent of the choice of $f$ and
gives a $\fg\otimes\C((t))\oplus\C\bk$-module structure on $Z^\infty$
with $\bk$ acting as scalar $\ell$. Restricted to Lie subalgebra
$\hat\fg$ of $\fg\otimes\C((t))\oplus\C\bk$, we have on $Z^\infty$
a structure of $\hat\fg$-module of level $\ell$.
We will denote this
$\hat\fg$-module by $\circ_C(W_1,W_2,\cdots,W_m)$, and when $m=2$,
simply by $W_1\circ_C W_2$.

Finally we need:
\begin{defi}{\rm
Given a $\hat\fg$-module $W$, consider $\hom(W,\C)$ as a
$\hat\fg$-module with the actions given by
\begin{equation}\label{W*2}
((a\otimes t^n)\lambda)(w)=-\lambda((a\otimes(-t)^{-n})w),
\quad (\bk\cdot\lambda)(w)=\lambda(\bk w)
\end{equation}
for $v\in{\mathfrak a}$, $\lambda \in W^*$ and $w\in W$. The {\em
contragredient module $D(W)$} of $W$ is defined by
\begin{eqnarray*}
D(W)&=&\{\lambda\in\hom(W,\C)\;|\;\mbox{there is }N\in\N\mbox{ such
that for any }a_1, \dots, a_N\in \fg,\\
&&\hspace{1em}(a_1\otimes t)\cdots(a_N\otimes t)\lambda=0\},
\end{eqnarray*}
a $\hat\fg$-submodule of $\hom(W,\C)$.  }
\end{defi}

The $\hat\fg$-module $D(\circ_C(W_1,W_2,\cdots,W_m))$ is defined to be
the desired tensor product of $W_1, \dots, W_m$.

Now we recall the category ${\cal O}_\kappa$ for a complex number
$\kappa\notin\Q_{\geq 0}$ from \cite{KL2}.

Let $M$ be a module for $\hat\fg_{(\leq 0)}$ satisfying the condition
that $\bk$ acts as a scalar $\ell$. Then the induced $\hat\fg$-module
\begin{equation}\label{Nkappa}
U(\hat\fg)\otimes_{U(\hat\fg_{(\leq 0)})}M
\end{equation}
is a $\hat\fg$-module of level $\ell$. In case $M$ is
finite-dimensional, restricted, and the subalgebra $\fg_{(-)}$ acts
nilpotently, the corresponding induced module is called a {\em
generalized Weyl module}.

\begin{defi}{\rm
Given complex number $\kappa\notin\Q_{\geq 0}$, the category ${\cal
O}_\kappa$ is defined to be the full subcategory of the category of
$\hat\fg$-modules whose objects are quotients of some generalized Weyl
module of level $\ell=\kappa-h$, where $h$ is the dual Coxeter number
of $\fg$. }
\end{defi}

For an object $W$ of ${\cal O}_\kappa$, define the {\em Segal-Sugawara
operator} $L_k: W\to W$ by
\begin{equation}\label{Lkdef}
L_k(w)=\frac{1}{2\kappa}\sum_{j\geq -k/2}\sum_p(c_p\otimes t^{-j})
(c_p\otimes t^{j+k})w+\frac{1}{2\kappa}\sum_{j<-k/2}\sum_p
(c_p\otimes t^{j+k})(c_p\otimes t^{-j})w
\end{equation}
where $\{c_p\}$ is an orthonormal basis of $\fg$. For any $n\in\C$
denote by $W_{[n]}$ the generalized eigenspace for $L_0$ with respect
to the eigenvalue $n$. Then it is shown in \cite{KL2} that
$W=\coprod_{n\in\C}W_{[n]}$ and $\dim W_{[n]}<\infty$.

The following result is proved in \cite{KL2}
\begin{theo}\label{closethm}
For $W_1$ and $W_2$ in ${\cal O}_\kappa$, $W_1\circ_{C} W_2$ is
again an object of ${\cal O}_\kappa$; the functor $D(\cdot)$ is closed
in ${\cal O}_\kappa$, and furthermore, as a vector space
\[
D(W)=\coprod_{n\in\C}(W_{[n]})^*.
\]
\end{theo}

Let $z$ be a nonzero complex
number. Consider the Riemann sphere $C$ with punctures $p_0=z$,
$p_1=\infty$ and $p_2=0$ and local coordinates given by
$\varphi_0(\epsilon)=\epsilon-z$, $\varphi_1(\epsilon)=1/\epsilon$ and
$\varphi_2(\epsilon)=\epsilon$ at $p_0$, $p_1$ and $p_2$,
respectively.  A Riemann sphere equipped with these
punctures and local coordinates is denoted by $Q(z)$ as in
\cite{tensor1}.

By (\ref{hom2}) and (\ref{aaawwwc}), the action of $\Gamma_R$ on
$W_1\otimes W_2$ associated with $Q(z)$ is given by $\bk$ acting as
$-\ell$ and
\begin{eqnarray}
(a\otimes f)(w_{(1)}\otimes
w_{(2)})&=&(a\otimes\iota_{\infty}f)(w_{(1)})\otimes w_{(2)}\nno\\
&&+w_{(1)}\otimes(a\otimes\iota_{0}f)(w_{(2)})
\end{eqnarray}
for $a\in\fg$, $f\in R$, $w_{(1)}\in W_1$ and $w_{(2)}\in W_2$. Hence
by (\ref{Liedual}) the action of $\Gamma_R$ on $(W_1\otimes W_2)^*$ is
given by $\bk$ acting as $\ell$ and
\begin{eqnarray*}
\lefteqn{\langle(a\otimes f)(\lambda),w_{(1)}\otimes w_{(2)}\rangle
=-\langle\lambda,(a\otimes f)(w_{(1)}\otimes w_{(2)})\rangle}\\
&&\hspace{1em}=-\langle\lambda,(a\otimes\iota_{\infty}f)
(w_{(1)})\otimes w_{(2)}+
w_{(1)}\otimes(a\otimes\iota_{0}f)(w_{(2)})\rangle
\end{eqnarray*}
for $a\in\fg$, $f\in R$, $\lambda\in(W_1\otimes W_2)^*$, $w_{(1)}\in
W_1$ and $w_{(2)}\in W_2$. This gives an action of the Lie algebra
$\fg\otimes\C[t,t^{-1},(z+t)^{-1}]\oplus\C\bk$.
(Note that $\C[t,t^{-1},(z+t)^{-1}]=\iota_{z}R$.)
In
particular, for $f=t^n$, $n\in\Z$, we have $\iota_{z}t^n=(z+t)^n$, and
\begin{eqnarray}
\lefteqn{\langle(a\otimes(z+t)^n)(\lambda),w_{(1)}\otimes
w_{(2)}\rangle}\nno\\
&&\hspace{1em}=-\langle\lambda,(a\otimes\iota_{\infty}t^n)(w_{(1)})
\otimes w_{(2)}+w_{(1)}\otimes(a\otimes\iota_{0}t^n)(w_{(2)})
\rangle\nno\\
&&\hspace{1em}=-\langle\lambda,(a\otimes t^{-n})(w_{(1)})\otimes
w_{(2)}+w_{(1)}\otimes(a\otimes t^n)(w_{(2)})\rangle
\end{eqnarray}
for any $a\in\fg$, $\lambda\in(W_1\otimes W_2)^*$,
$w_{(1)}\in W_1$ and $w_{(2)}\in W_2$; and in case $f=(t-z)^n$, $n\in\Z$,
we have $\iota_{z}(t-z)^n=t^n$, and
\begin{eqnarray}
\lefteqn{\langle(a\otimes t^n)(\lambda),w_{(1)}\otimes
w_{(2)}\rangle}\nno\\
&&\hspace{1em}=-\langle\lambda,(a\otimes\iota_{\infty}(t-z)^n)
(w_{(1)})\otimes w_{(2)}+w_{(1)}\otimes(a\otimes\iota_{0}(t-z)^n)
(w_{(2)})\rangle\nno\\
&&\hspace{1em}=-\langle\lambda,(a\otimes\iota_+(t^{-1}-z)^n)
(w_{(1)})\otimes w_{(2)}+w_{(1)}\otimes(a\otimes\iota_+(t-z)^n)
(w_{(2)})\rangle\nno\\
&&\hspace{1em}=-\Big\langle\lambda,\sum_{i\in\N}{n\choose i}(-z)^i
(a\otimes t^{i-n})(w_{(1)})\otimes w_{(2)}+\nno\\
&&\hspace{3em}+w_{(1)}\otimes\sum_{i\in\N}{n\choose i}(-z)^{n-i}
(a\otimes t^i)(w_{(2)})\Big\rangle
\end{eqnarray}
for any $a\in\fg$, $\lambda\in(W_1\otimes W_2)^*$,
$w_{(1)}\in W_1$ and $w_{(2)}\in W_2$.

By the general construction we now have
\begin{eqnarray}\label{circ_Q}
W_1\circ_{Q(z)}W_2&=&\{\lambda\in(W_1\otimes W_2)^*\;|\; \mbox{for
some }N\in\N,\; \xi_1\xi_2\cdots\xi_N\lambda=0\nno\\
&&\hspace{1em}\mbox{ for any }\xi_1,\xi_2,\dots,\xi_N\in\fg\otimes
t\C[t,(z+t)^{-1}] \}.
\end{eqnarray}
The tensor product is then defined as
\[
D(W_1\circ_{Q(z)}W_2).
\]

\setcounter{equation}{0} \setcounter{rema}{0}
\section{Tensor product for modules for a conformal vertex algebra}

In this section we recall the construction and some results in
the tensor product theory for suitable module categories for a
conformal vertex algebra from \cite{tensor1}--\cite{tensor3},
\cite{tensor4} and \cite{HLZ1}, \cite{HLZ2}.

We assume the reader is familiar with the material in \cite{FLM} and
\cite{FHL}, such as the language of formal calculus, the notion of
vertex operator algebra and their modules, etc.  Results from
\cite{tensor1}--\cite{tensor3}, \cite{tensor4} and \cite{HLZ1},
\cite{HLZ2} will be recalled without proof. We refer the reader to
these papers for details.

We will focus on the ``$Q(z)$-tensor product'' in this section, due to
the fact that tensor product constructed in \cite{KL1} corresponds to
$Q(1)$.

In \cite{HLZ1} and \cite{HLZ2}, for an abelian group 
$A$ and an abelian group 
$\tilde{A}$ containing $A$ as a subgroup, the notions of strongly 
$A$-graded conformal vertex
algebra and strongly $\tilde{A}$-graded generalized modules 
were introduced.  The vertex operator algebras and their (ordinary) modules are
exactly the conformal vertex algebras and their (ordinary) modules
that are strongly graded with respect to the trivial group (see 
Remark 2.4 in \cite{HLZ1} and Remarks 2.24 and 2.27 in \cite{HLZ2}). 
In the present paper, we shall work 
in the special case of \cite{HLZ1} and \cite{HLZ2} in which 
$A$ and $\tilde{A}$ are trivial. 

Let $V$ be a vertex operator algebra, that is, as a special case 
of the general theory developed in \cite{HLZ1} and \cite{HLZ2},
a strongly $A$-graded conformal vertex algebra with trivial $A$. 
A {\it generalized $V$-module} is a strongly $\tilde{A}$-graded 
generalized module in the sense of \cite{HLZ1} and \cite{HLZ2}
with trivial $\tilde{A}$.
It can also be defined directly in the same way as a 
$V$-module except that instead of
being the direct sum of eigenspaces for the operator $L(0)$, it is
assumed to be the direct sum of {\em generalized} eigenspaces for
$L(0)$.

Recall from Definition 2.5 in  \cite{HLZ1} or
Definition 3.32 in \cite{HLZ2} 
that given a generalized $V$-module $(W,Y_W)$ with
\[
W=\coprod_{n\in\C}W_{[n]}
\]
where $W_{[n]}$ is the generalized eigenspace for $L(0)$ with respect
to the eigenvalue $n$, its {\em contragredient module} is the vector
space
\[
W'=\coprod_{n\in\C}(W_{[n]})^*
\]
equipped with the vertex operator map $Y'$ defined by
\[
\langle Y'(v,x)w',w\rangle=\langle w',Y^o_W(v,x)w\rangle,
\]
for any $v\in V$, $w'\in W'$ and $w\in W$, where
\[
Y^o_W(v,x)=Y_W(e^{xL(1)}(-x^{-2})^{L(0)}v,x^{-1}),
\]
for any $v\in V$, is the {\em opposite vertex operator} (cf. \cite{FHL}). We will use
the standard notation
\[
\overline W=\prod_{n\in\C}W_{[n]},
\]
the formal completion of $W$ with respect to the $\C$-grading.

Fix a nonzero complex number $z$. The concept of $Q(z)$-intertwining
map is defined as follows (see Definition 4.17 of \cite{HLZ1} or
Definition 4.32 of \cite{HLZ2}):
\begin{defi}{\rm
Let $W_1$, $W_2$ and $W_3$ be generalized modules for a vertex
operator algebra $V$.  A {\it $Q(z)$-intertwining map of type
${W_3\choose W_1\,W_2}$} is a linear map $I: W_1\otimes W_2 \to
\overline{W}_3$ such that the following conditions are satisfied: the
{\em lower truncation condition:} for any elements $w_{(1)}\in W_1$,
$w_{(2)}\in W_2$, and any $n\in {\mathbb C}$,
\begin{equation}\label{imq:ltc}
\pi_{n-m}I(w_{(1)}\otimes w_{(2)})=0\;\;\mbox{ for }\;m\in {\mathbb
N}\;\mbox{ sufficiently large;}
\end{equation}
and the {\em Jacobi identity}:
\begin{eqnarray}\label{imq:def}
\lefteqn{z^{-1}\delta\left(\frac{x_1-x_0}{z}\right)
Y^o_3(v, x_0)I(w_{(1)}\otimes w_{(2)})=}\nonumber\\
&&=x_0^{-1}\delta\left(\frac{x_1-z}{x_0}\right)
I(Y_1^{o}(v, x_1)w_{(1)}\otimes w_{(2)})\nonumber\\
&&\hspace{2em}-x_0^{-1}\delta\left(\frac{z-x_1}{-x_0}\right)
I(w_{(1)}\otimes Y_2(v, x_1)w_{(2)})
\end{eqnarray}
for $v\in V$, $w_{(1)}\in W_1$ and $w_{(2)}\in W_2$ (note that the
left-hand side of (\ref{imq:def}) is meaningful because any infinite
linear combination of $v_n$ of the form $\sum_{n<N}a_nv_n$ ($a_n\in
{\mathbb C}$) acts on any $I(w_{(1)}\otimes w_{(2)})$, due to
(\ref{imq:ltc})). }
\end{defi}

Given generalized $V$-modules $W_1$ and $W_2$, we first have the
notion of a $Q(z)$-product, as follows (see Definition 4.39 of \cite{HLZ2}):
\begin{defi}{\rm
Let $W_1$ and $W_2$ be generalized $V$-modules.  A {\it
$Q(z)$-product of $W_1$ and $W_2$} is a generalized $V$-module $(W_3,
Y_3)$ together with a $Q(z)$-intertwining map $I_3$ of type
${W_3}\choose {W_1W_2}$, We denote it by $(W_3, Y_3; I_3)$ or simply
by $(W_3, I_3)$.  Let $(W_4,Y_4; I_4)$ be another $Q(z)$-product of
$W_1$ and $W_2$. A {\em morphism} from $(W_3, Y_3; I_3)$ to $(W_4,
Y_4; I_4)$ is a module map $\eta$ from $W_3$ to $W_4$ such that
\[
I_4=\overline{\eta}\circ I_3.
\]
where $\overline{\eta}$ is the natural map {from} $\overline{W}_3$ to
$\overline{W}_4$ uniquely extending $\eta$. }
\end{defi}

Let ${\cal C}$ be a full subcategory of the category of generalized
$V$-modules. The notion of $Q(z)$-tensor product of $W_1$ and $W_2$ in
${\cal C}$ is defined in term of a universal property as follows (see
Definition 4.40 of \cite{HLZ2}):
\begin{defi}\label{Qtpdef}{\rm
For $W_1, W_2\in \ob{\cal C}$, a {\em $Q(z)$-tensor product of $W_1$
and $W_2$ in ${\cal C}$} is a $Q(z)$-product $(W_0, Y_0; I_0)$ with
$W_0\in{\rm ob\,}{\cal C}$ such that for any $Q(z)$-product $(W,Y;I)$
with $W\in{\rm ob\,}{\cal C}$, there is a unique morphism from $(W_0,
Y_0; I_0)$ to $(W,Y;I)$. Clearly, a $Q(z)$-tensor product of $W_1$ and
$W_2$ in ${\cal C}$, if it exists, is unique up to a unique
isomorphism. In this case we will denote it as $(W_1\boxtimes_{Q(z)}
W_2, Y_{Q(z)}; \boxtimes_{Q(z)})$ and call the object
$(W_1\boxtimes_{Q(z)} W_2, Y_{Q(z)})$ the {\em $Q(z)$-tensor product
module of $W_1$ and $W_2$ in ${\cal C}$}. We will skip the term ``in
${\cal C}$'' if the category ${\cal C}$ under consideration is clear
in context. }
\end{defi}

Now we recall the construction of the $Q(z)$-tensor product from
Section 5.3 of \cite{HLZ2}, which generalizes that in
\cite{tensor1}--\cite{tensor3}.  Let $W_1$ and $W_2$ be generalized
$V$-modules. We first have the following linear action $\tau_{Q(z)}$
of the space
\[
V\otimes \iota_{+}\C[t, t^{-1}, (z+t)^{-1}]
\]
on $(W_1\otimes W_2)^*$:
\begin{eqnarray}\label{tauQ}
\lefteqn{\left(\tau_{Q(z)}
\left(z^{-1}\delta\left(\frac{x_1-x_0}{z}\right) Y_{t}(v,
x_0)\right)\lambda\right)(w_{(1)}\otimes w_{(2)})}\nno\\
&&=x^{-1}_0\delta\left(\frac{x_1-z}{x_0}\right) \lambda(Y_1^{o}(v,
x_1)w_{(1)}\otimes w_{(2)})\nno\\
&&\hspace{2em}-x_0^{-1}\delta\left(\frac{z-x_1}{-x_0}\right)
\lambda(w_{(1)}\otimes Y_2(v, x_1)w_{(2)}).
\end{eqnarray}
for $v\in V$, $\lambda\in (W_1\otimes W_2)^{*}$, $w_{(1)}\in W_1$,
$w_{(2)}\in W_2$, where
\begin{equation}\label{Y_t}
Y_t(v,x)=v\otimes t^{-1}\delta\bigg(\frac{t}{x}\bigg).
\end{equation}
This includes an action $Y'_{Q(z)}$ of $V\otimes\C[t,t^{-1}]$ defined
by
\begin{equation}\label{Y'Q}
Y'_{Q(z)}(v,x)=\tau_{Q(z)}(v\otimes t^{-1}\delta\bigg(
\frac{t}{x}\bigg)),
\end{equation}
that is, by taking $\res_{x_1}$ in (\ref{tauQ}),
\begin{eqnarray}
\lefteqn{(Y'_{Q(z)}(v,x_0)\lambda)(w_{(1)} \otimes w_{(2)})=}\nno\\
&&= \lambda(Y^o_1(v,x_0 + z)w_{(1)} \otimes w_{(2)}) \nno\\
&&\hspace{2em}- \res_{x_1}x^{-1}_0 \delta
\left(\frac{z-x_1}{-x_0}\right)\lambda(w_{(1)} \otimes
Y_2(v,x_1)w_{(2)}).
\end{eqnarray}
We also have the operators $L'_{Q(z)}(n)$, $n\in\Z$ defined by
\[
Y'_{Q(z)}(\omega,x)=\sum_{n\in\Z}L'_{Q(z)}(n)x^{-n-2}.
\]

We have the following construction of $W_1\hboxtr_{Q(z)}W_2$, a
subspace of $(W_1\otimes W_2)^*$ (see Definition 5.60 and Theorem 5.74
of \cite{HLZ2}):
\begin{defi}{\rm
Given $W_1$ and $W_2$ as above, the vector space
$W_1\hboxtr_{Q(z)}W_2$ consists of all the elements
$\lambda\in(W_1\otimes W_2)^*$ satisfying the following two conditions:

\begin{description}
\item{\bf The $Q(z)$-compatibility condition}

(a) The {\em lower truncation condition}: For all $v\in V$, the formal
Laurent series $Y'_{Q(z)}(v, x)\lambda$ involves only finitely many
negative powers of $x$.

(b) The following formula holds:
\begin{eqnarray}\label{comp}
\lefteqn{\tau_{Q(z)}\left(z^{-1}\delta\left(\frac{x_1-x_0}{z}\right)
Y_{t}(v, x_0)\right)\lambda=}\nno\\
&&=z^{-1}\delta\left(\frac{x_1-x_0}{z}\right)
Y'_{Q(z)}(v, x_0)\lambda  \;\;\mbox{ for all }\;v\in V.
\end{eqnarray}

\item{\bf The $Q(z)$-local grading-restriction condition}

(a) The {\em grading condition}: $\lambda$ is a (finite) sum of
generalized eigenvectors of $(W_1\otimes W_2)^{*}$ for the operator
$L'_{Q(z)}(0)$.

(b) The smallest subspace $W_{\lambda}$ of $(W_1\otimes W_2)^*$
containing $\lambda$ and stable under the component operators
$\tau_{Q(z)}(v\otimes t^{n})$ of the operators $Y'_{Q(z)}(v,x)$ for
$v\in V$, $n\in {\mathbb Z}$, have the properties:
\begin{eqnarray}
&\dim(W_\lambda)_{[n]}<\infty&\\
&(W_\lambda)_{[n+k]}=0\;\;\mbox{ for}\;k\in {\mathbb Z}
\;\mbox{ sufficiently negative};&
\end{eqnarray}
for any $n\in {\mathbb C}$, where the subscripts denote the ${\mathbb
C}$-grading by $L'_{Q(z)}(0)$-eigenvalues.
\end{description}
}
\end{defi}

The importance of the space $W_1\hboxtr_{Q(z)}W_2$ is given by the
following theorem from \cite{tensor1} and its generalization in
\cite{HLZ2} (see Theorems 5.70, 5.71, 5.72, 5.73 and 5.74 of 
\cite{HLZ2}): 
\begin{theo} The vector space $W_1\hboxtr_{Q(z)}W_2$ is closed
under the action $Y'_{Q(z)}$ of $V$ and the Jacobi identity holds on
$W_1\hboxtr_{Q(z)}W_2$. Furthermore, the $Q(z)$-tensor product of
$W_1$ and $W_2$ in ${\cal C}$ exists if and only if
$W_1\hboxtr_{Q(z)}W_2$ equipped with $Y'_{Q(z)}$ is an object of
${\cal C}$, and in this case, this $Q(z)$-tensor product is the
contragredient module of $(W_1\hboxtr_{Q(z)}W_2, Y'_{Q(z)})$.
\end{theo}

\section{Strong lower truncation condition}
In this section we define what we shall call the strong lower truncation condition and
prove its equivalence to the compatibility condition.

By using (\ref{Y_t}), (\ref{Y'Q}) and the delta function identity
\[
z^{-1}\delta\left(\frac{x_1-x_0}{z}\right)=
x_1^{-1}\delta\left(\frac{z+x_0}{x_1}\right),
\]
formula (\ref{comp}) in the compatibility condition can be written as
\begin{eqnarray*}
\lefteqn{\tau_{Q(z)}\left(x_1^{-1}\delta\left(\frac{z+x_0}{x_1}\right)
Y_{t}(v,x_0)\right)\lambda=}\nno\\
&&=x_1^{-1}\delta\left(\frac{z+x_0}{x_1}\right)\tau_{Q(z)}(Y_{t}(v,x_0)
)\lambda.
\end{eqnarray*}
Taking the coefficient of $x_1^n$ for $n\in\Z$ of both sides we get
\[
\tau_{Q(z)}((z+x_0)^{-n-1}Y_{t}(v,x_0))\lambda=
(z+x_0)^{-n-1}Y'_{Q(z)}(v,x_0)\lambda.
\]
or, by using (\ref{Y_t}) and the property of the $\delta$-function
\[
t^{-1}\delta\bigg(\frac{t}{x_0}\bigg)f(x_0)=t^{-1}\delta\bigg(\frac{t}{x_0}\bigg)f(t)
\]
for formal series $f(x_0)$ we have
\begin{equation}\label{compn}
\tau_{Q(z)}((z+t)^{-n-1}Y_{t}(v,x_0))\lambda=
(z+x_0)^{-n-1}Y'_{Q(z)}(v,x_0)\lambda.
\end{equation}
Further taking the coefficient of $x_0^{-m-1}$ for $m\in\Z$ this becomes
\begin{equation}\label{comp1}
\tau_{Q(z)}(v\otimes (z+t)^{-n-1}t^m)\lambda=\sum_{i\in\N}{-n-1
\choose i}z^{-n-1-i}\tau_{Q(z)}(v\otimes t^{m+i})\lambda.
\end{equation}

We have:
\begin{propo}\label{obser}
Let $W_1$ and $W_2$ be modules for $V$ as a vertex algebra,
$\lambda\in (W_1\otimes W_2)^*$, $v\in V$ and $n\in\N$. Then
\[
\tau_{Q(z)}((z+t)^{-n-1} Y_{t}(v,x_0))\lambda
\]
is lower truncated in $x_0$ if and only if $Y'_{Q(z)}(v, x_0)\lambda$
is lower truncated in $x_0$ and (\ref{compn}) holds.
\end{propo}
\pf The ``if'' part is obvious. For the ``only if'' part,
suppose that 
\[
\tau_{Q(z)}((z+t)^{-n-1} Y_{t}(v,x_0))\lambda=
\tau_{Q(z)}((z+x_0)^{-n-1} Y_{t}(v,x_0))\lambda
\]
is lower truncated in $x_0$, then so is
\begin{eqnarray}\label{yy}
&&(z+x_0)^{n+1}\tau_{Q(z)}((z+x_0)^{-n-1} Y_{t}(v,x_0))\lambda\nno\\
&=&\tau_{Q(z)}((z+x_0)^{n+1}(z+x_0)^{-n-1} Y_{t}(v,x_0))\lambda\nno\\
&=&\tau_{Q(z)}(Y_{t}(v,x_0))\lambda\nno\\
&=&Y'_{Q(z)}(v,x_0)\lambda.
\end{eqnarray}
(Note that here we need the existence of the triple product
\[
(z+x_0)^{n+1}(z+x_0)^{-n-1} Y_{t}(v,x_0)
\]
which can be seen by, for example, observing that the coefficient of
each power of $t$ exists.) That is, $Y'_{Q(z)}(v,x_0)\lambda$ is lower
truncated in $x_0$.

By (\ref{yy}) we also have
\begin{eqnarray}\label{=0}
&&(z+x_0)^{n+1}(\tau_{Q(z)}((z+x_0)^{-n-1} Y_{t}(v,x_0))\lambda-
(z+x_0)^{-n-1}Y'_{Q(z)}(v,x_0)\lambda)\nno\\
&=&(z+x_0)^{n+1}\tau_{Q(z)}((z+x_0)^{-n-1} Y_{t}(v,x_0))\lambda-
Y'_{Q(z)}(v,x_0)\lambda\nno\\
&=&Y'_{Q(z)}(v,x_0)\lambda-Y'_{Q(z)}(v,x_0)\lambda\nno\\
&=&0.
\end{eqnarray}
Since both terms of
\begin{equation}\label{0}
\tau_{Q(z)}((z+x_0)^{-n-1} Y_{t}(v,x_0))\lambda-
(z+x_0)^{-n-1}Y'_{Q(z)}(v,x_0)\lambda
\end{equation}
are lower truncated in $x_0$, (\ref{=0}) implies that (\ref{0}) is
equal to $0$, as desired.
\epfv

Now we define the strong lower truncation condition:
\begin{defi}{\rm
Let $v\in V$. An element $\lambda$ in $(W_1\otimes W_2)^*$ is said to
satisfy the {\em strong lower truncation condition with respect to
$v\in V$} if there exists $N\in\N$ depending on $v$ and $\lambda$
such that
\begin{equation}\label{compeqi}
(\tau_{Q(z)}(v\otimes t^m(z+t)^{-n-1}))\lambda=0
\end{equation}
for all $m\geq N$ and $n\in\Z$.
We say that $\lambda$ satisfies the {\em strong
lower truncation condition} if it satisfies the strong lower
truncation condition with respect to every vector in $V$.  }
\end{defi}

As a consequence of Proposition \ref{obser} we have the following
equivalent condition for the $Q(x)$-compatibility condition:
\begin{propo}\label{coro}
Let $\lambda\in (W_1\otimes W_2)^*$. Then $\lambda$ satisfies the
$Q(z)$-compatibility condition if
and only if $\lambda$ satisfies the strong lower truncation condition.
\end{propo}
\pf Suppose that $\lambda$ satisfies 
the $Q(z)$-compatibility condition. Then (\ref{comp1})
holds for any $v\in V$ and $m,n\in\Z$. But then by part (a) of the
compatibility condition we see that the right-hand side of
(\ref{comp1}) is $0$ when $m$ is large enough, independent of
$n$. This proves half of the statement. The other half follows
directly from Proposition \ref{obser}. \epfv

We will need:
\begin{lemma}\label{lemma:slt}
Let $N\in\N$. For $v\in V$ and $\lambda\in(W_1\otimes W_2)^*$, 
(\ref{compeqi}) holds
for any $m\geq N$ and $n\in\Z$ if and only if for any $w_{(1)}\in W_1$
and $w_{(2)}\in W_2$,
\[
(x_1-z)^N\lambda(Y_1^{o}(v, x_1)w_{(1)}\otimes w_{(2)}-w_{(1)}\otimes
Y_2(v, x_1)w_{(2)})=0.
\]
\end{lemma}
\pf By definition, (\ref{compeqi}) holds for any $m\geq N$ and $n\in\Z$
if and only if all powers of $x_0$ with nonzero coefficients in
\[
\tau_{Q(z)}\left(z^{-1}\delta\left(\frac{x_1-x_0}{z}\right) Y_{t}(v,
x_0)\right)\lambda
\]
are at least $-N$. By definition of the action $\tau_{Q(z)}$ in
(\ref{tauQ}) this in turn is equivalent to the condition that
$\res_{x_0}x_0^m$ of the right-hand side of (\ref{tauQ}) is $0$ for
any $m\geq N$, $w_{(1)}\in W_1$ and $w_{(2)}\in W_2$. The statement
now follows from this.
\epfv

The following result in formal calculus will be handy for us:
\begin{lemma}\label{lem:==}
Let $x$ and $y$ be formal variables, $\xi$ be either a formal variable
or a complex number. Let $K\in\N$ and let $f_k(x,\xi)$, $k\in\Z$ be a
sequence of formal series with coefficients in some vector space
satisfying the condition that $f_k(x,\xi)=0$ for any $k\geq
K$. Suppose that for $N_1, N_2\in\N$ we have
\begin{equation}\label{++}
(x+\xi)^{N_1}(x+y+\xi)^{N_2}\sum_{n\in\Z}f_n(x,\xi)y^{-n-1}=0.
\end{equation}
Then
\begin{equation}\label{++=}
(x+\xi)^{N_1+N_2+s}f_{K-1-s}(x,\xi)=0
\end{equation}
for any $s\in\N$.
\end{lemma}
\pf By taking coefficient of powers $y$, we see that (\ref{++}) is
equivalent to
\begin{equation}\label{++.}
\sum_{i=0}^{N_2}{N_2\choose i}(x+\xi)^{N_1+N_2-i}f_{i+k}(x,\xi)=0
\end{equation}
for any $k\in\Z$. Since $f_k(x,\xi)=0$ for any $k\geq K$, setting
$k=K-1$ in (\ref{++.})  we obtain
$(x+\xi)^{N_1+N_2}f_{K-1}(x,\xi)=0$. Now (\ref{++=}) follows by
induction on $s$, as follows: If the case $0, 1, \dots, s-1$ holds,
then by setting
$k=K-1-s$ in (\ref{++.}) and multiplying both sides by $(x+\xi)^s$ we
see that only the $i=0$ term remains and must equal the right-hand
side $0$, that is, $(x+\xi)^{N_1+N_2+s}f_{K-1-s}(x,\xi)=0$ also holds.
\epfv

We now have:
\begin{propo}\label{p:gene}
Let $u, v\in V$ and $\lambda\in(W_1\otimes W_2)^*$. Suppose that
$\lambda$ satisfies the strong lower truncation condition with respect
to both $u$ and $v$. Then for any $k\in\Z$, $\lambda$ satisfies the
strong lower truncation condition with respect to $u_k v$. More
precisely, let $N_1$ be the integer such that (\ref{compeqi}) holds
with $v$ replaced by $u$ for any $m\geq N_1$ and $n\in\Z$ and let
$N_2$ be the corresponding number for $v$, then (\ref{compeqi}) holds
with $v$ replaced by $u_k v$ for any $m\geq N_1+N_2+K-1-k$ and
$n\in\Z$, where $K\in\N$ is such that $u_n v=0$ for any $n\geq K$.
\end{propo}
\pf By assumption and Lemma \ref{lemma:slt} we have
\[
(x_1-z)^{N_1}\lambda(Y_1^{o}(u, x_1)w_{(1)}\otimes w_{(2)}-w_{(1)}
\otimes Y_2(u, x_1)w_{(2)})=0.
\]
and
\[
(x_1-z)^{N_2}\lambda(Y_1^{o}(v, x_1)w_{(1)}\otimes w_{(2)}-w_{(1)}
\otimes Y_2(v, x_1)w_{(2)})=0.
\]
for any $w_{(1)}\in W_1$ and $w_{(2)}\in W_2$. For $k\in\Z$, we need a formula
similar to either of these, with $u$ or $v$ replaced by $u_k v$. We
derive as follows: First, for formal variables $x-1$, $y_0$ and $y_1$,
using the above identities we have:
\begin{eqnarray*}
\lefteqn{(x_1-z)^{N_2}(x_1+y_0-z)^{N_1}\lambda(y^{-1}_0\delta\left({y_1
-x_1\over y_0}\right)Y_1^o(v,x_1)Y_1^o(u,y_1)w_{(1)}\otimes w_{(2)})}\\
&&=y^{-1}_0\delta\left({y_1-x_1\over y_0}\right)(x_1-z)^{N_2}
(x_1+y_0-z)^{N_1}\lambda(Y_1^o(v,x_1)Y_1^o(u,y_1)w_{(1)}\otimes w_{(2)})\\
&&=y^{-1}_0\delta\left({y_1-x_1\over y_0}\right)(x_1-z)^{N_2}
(x_1+y_0-z)^{N_1}\lambda(Y_1^o(u,y_1)w_{(1)}\otimes Y_2(v,x_1)w_{(2)})\\
&&=y^{-1}_0\delta\left({y_1-x_1\over y_0}\right)(x_1-z)^{N_2}
(y_1-z)^{N_1}\lambda(Y_1^o(u,y_1)w_{(1)}\otimes Y_2(v,x_1)w_{(2)})\\
&&=y^{-1}_0\delta\left({y_1-x_1\over y_0}\right)(x_1-z)^{N_2}
(y_1-z)^{N_1}\lambda(w_{(1)}\otimes Y_2(u,y_1)Y_2(v,x_1)w_{(2)})\\
&&=y^{-1}_0\delta\left({y_1-x_1\over y_0}\right)(x_1-z)^{N_2}
(x_1+y_0-z)^{N_1}\lambda(w_{(1)}\otimes Y_2(u,y_1)Y_2(v,x_1)w_{(2)})\\
&&=(x_1-z)^{N_2}(x_1+y_0-z)^{N_1}\lambda(w_{(1)}\otimes y^{-1}_0
\delta\left({y_1-x_1\over y_0}\right)Y_2(u,y_1)Y_2(v,x_1)w_{(2)}).
\end{eqnarray*}
On the other hand,
\begin{eqnarray*}
\lefteqn{(x_1-z)^{N_2}(x_1+y_0-z)^{N_1}\lambda(y^{-1}_0\delta\left({x_1-
y_1\over -y_0}\right)Y_1^o(u,y_1)Y^o_1(v,x_1)w_{(1)}\otimes w_{(2)})}\\
&&=y^{-1}_0\delta\left( {x_1-y_1\over -y_0}\right)(x_1-z)^{N_2}
(x_1+y_0-z)^{N_1}\lambda(Y_1^o(u,y_1)Y^o_1(v,x_1)w_{(1)}\otimes w_{(2)})\\
&&=y^{-1}_0\delta\left( {x_1-y_1\over -y_0}\right)(x_1-z)^{N_2}
(y_1-z)^{N_1}\lambda(Y_1^o(u,y_1)Y^o_1(v,x_1)w_{(1)}\otimes w_{(2)})\\
&&=y^{-1}_0\delta\left( {x_1-y_1\over -y_0}\right)(x_1-z)^{N_2}
(y_1-z)^{N_1}\lambda(Y^o_1(v,x_1)w_{(1)}\otimes Y_2(u,y_1)w_{(2)})\\
&&=y^{-1}_0\delta\left( {x_1-y_1\over -y_0}\right)(x_1-z)^{N_2}
(x_1+y_0-z)^{N_1}\lambda(Y^o_1(v,x_1)w_{(1)}\otimes Y_2(u,y_1)w_{(2)})\\
&&=y^{-1}_0\delta \left({x_1-y_1\over -y_0}\right)(x_1-z)^{N_2}
(x_1+y_0-z)^{N_1}\lambda(w_{(1)}\otimes Y_2(v,x_1)Y_2(u,y_1)w_{(2)})\\
&&=(x_1-z)^{N_2}(x_1+y_0-z)^{N_1}\lambda(w_{(1)}\otimes y^{-1}_0
\delta\left({x_1-y_1\over -y_0}\right) Y_2(v,x_1)Y_2(u,y_1)w_{(2)}).
\end{eqnarray*}
Taking difference of these two equalities, using the Jacobi identity
and the opposite Jacobi identity we have
\begin{eqnarray*}
\lefteqn{(x_1-z)^{N_2}(x_1+y_0-z)^{N_1}\lambda(x^{-1}_1\delta\left({y_1-
y_0\over x_1}\right)Y_1^o(Y(u,y_0)v,x_1)w_{(1)}\otimes w_{(2)})}\\
&&=(x_1-z)^{N_2}(x_1+y_0-z)^{N_1}\lambda(w_{(1)}\otimes x^{-1}_1
\delta\left({y_1-y_0\over x_1}\right) Y_2(Y(u,y_0)v,x_1)w_{(2)}).
\end{eqnarray*}

Applying $\res_{y_1}$ we have
\begin{eqnarray*}
\lefteqn{(x_1-z)^{N_2}(x_1+y_0-z)^{N_1}\cdot}\\
&&\hspace{3em}\lambda(Y_1^o(Y(u,y_0)v,x_1)w_{(1)}\otimes w_{(2)}-
w_{(1)}\otimes Y_2(Y(u,y_0)v,x_1)w_{(2)})=0.
\end{eqnarray*}
Let $K$ be a number such that $u_k v=0$ for all $k\geq K$. Then
Lemma \ref{lem:==} applies with $\xi$ being $-z$ and $f_k(x,\xi)$
being $\lambda(Y_1^o(u_k v,x_1)w_{(1)}\otimes w_{(2)}-w_{(1)}\otimes
Y_2(u_k v,x_1)w_{(2)})$. As a result we have
\[
(x_1-z)^{N_2+N_2+s}\lambda(Y_1^o(u_{K-1-s}v,x_1)w_{(1)}\otimes w_{(2)}-
w_{(1)}\otimes Y_2(u_{K-1-s}v,x_1)w_{(2)})=0,
\]
for any $s\in\N$. This is exactly what we need.
\epf

Combining Proposition \ref{coro} and Proposition \ref{p:gene} we
obtain:
\begin{theo}\label{th:main}
Let $V$ be a vertex algebra and $S$ a generating set for $V$ in the
sense that
\[
V={\rm span}\{a^{(k)}_{n_{k}}\cdots a^{(2)}_{n_2}a^{(1)}_{n_1}
a^{(0)}\,|\,a^{(0)},a^{(1)},\dots,a^{(k)}\in S,\;n_1,n_1,\dots,
n_k\in\Z, k\in\N\}.
\]
Let $W_1$ and $W_2$ be generalized $V$-modules. Then $\lambda\in
(W_1\otimes W_2)^*$ satisfies the $Q(z)$-compatibility condition if
and only if $\lambda$ satisfies the strong lower truncation condition
with respect to all elements of $S$. \epf
\end{theo}

\section{Generalized modules for vertex operator algebras associated
to $\hat\fg$}

In this section we recall the vertex operator algebra constructed from
suitable modules for the affine Lie algebra $\hat\fg$. Then we show
that the construction of the tensor products in Section 2 and Section 3
are equivalent, for modules in ${\cal O}_\kappa$. In particular, this
shows that the tensor product constructed by Kazhdan and Lusztig
satisfies the universal property in Definition \ref{Qtpdef}.  
Then we show that the objects of the category ${\cal
O}_\kappa$ are $C_{1}$-cofinite and quasi-finite dimensional.
These results together with a result 
in \cite{H3} imply that the conditions needed for applying 
the results obtained in \cite{HLZ1} and \cite{HLZ2}
are satisfied.  Hence we prove the existence of the associativity isomorphisms
and we obtain the braided tensor category structure.

Recall the complex semisimple Lie algebra $\fg$, the nondegenerate
invariant bilinear form $(\cdot, \cdot)$ on $\fg$ and the 
affine Lie algebra $\hat{\fg}$ in Section 2.

Given any $\fg$-module $U$ and any complex number $\ell$, consider $U$
as a $\hat\fg_{(\leq)}$-module with $\hat\fg_{(-)}$ acting trivially
and $\bk$ acting as the scalar $\ell$. Then the induced
$\hat\fg$-module
\[
{\rm Ind}^{\hat\fg}_\fg
(U)=U(\hat\fg)\otimes_{U(\hat\fg_{(\leq)})} U.
\]
is restricted and is of level $\ell$. When $U=\C$ is the trivial
$\fg$-module, we will write
\[
V_{\hat\fg}(\ell,0)={\rm Ind}^{\hat \fg}_\fg(\C).
\]

In particular, fix a Cartan subalgebra $\fh\subset \fg$ and a set of
positive roots, let $\fg=\fn_+\oplus\fh \oplus\fn_-$ be the
corresponding triangular decomposition. Let $U=L(\lambda)$ be the
irreducible highest weight $\fg$-module with highest weight
$\lambda\in\fh^*$. That is, $L(\lambda)$ is the quotient of
$V(\lambda)$ by its maximal proper submodule, where $V(\lambda)$ is
the $\fg$-module induced by the $(\fh\oplus \fn_+)$-module $\C
v_\lambda$ with
\[
hv_\lambda=\lambda(h)v_\lambda\;\mbox{ for }h\in\fh
\]
and $\fn_+$ acts on $v_\lambda$ as $0$. $L(\lambda)$ is finite
dimensional if and only if $\lambda$ is dominant integral in the sense
that
\[
\lambda(h_\alpha)=\frac{2(\lambda,\alpha)}{(\alpha,\alpha)}\in\N,
\;\mbox{ for }\alpha\in\Delta_+.
\]
In this case,
\[
L(\lambda)=U(\fn_-)/U(\fn_-)\alpha_i^{\lambda(h_{\alpha_i})+1}
\]
and we will denote the induced $\hat\fg$-module by $M(\ell,\lambda)$,
called the {\em Weyl module} for $\hat\fg$ with respect to $\lambda$.
Let $J(\ell,\lambda)$ be the maximal proper $\hat\fg$-submodule of
$M(\ell,\lambda)$. Then
$L(\ell,\lambda)=M(\ell,\lambda)/J(\ell,\lambda)$ is an irreducible
$\hat\fg$-module of level $\ell$.

We have (see \cite{KL2}):
\begin{propo}
(a) The operator $L(0)$ acts semisimply on $M(\ell,\lambda)$.

(b) The category ${\cal O}_\kappa$ consists of all the $\hat\fg$-modules
of level $\ell$ having a finite composition series all of whose
irreducible subquotients are of the form $L(\ell,\lambda)$ for various $\lambda$.
\end{propo}

Note that $V_{\hat\fg}(\ell,0)$ is spanned by the elements of the form
$a^{(1)}(-n_1)\cdots a^{(r)}(-n_r)1$, where $a^{(1)}, \dots,
a^{(r)}\in \fg$ and $n_1, \dots, n_r\in \Z_+$; here and below we use
$a(-n)$ to denote the representation image of $a\otimes t^{-n}$ for
$a\in \fg$ and $n\in \Z$.

The following theorem is well known:
\begin{theo} (\cite{FZ}; cf.\ \cite{LL})
There is a unique vertex algebra structure $(V_{\hat\fg}(\ell,0), Y, {\bf
1})$ on $V_{\hat\fg}(\ell,0)$ such that ${\bf 1}=1\in \C$ is the
vacuum vector and such that
\[
Y(a(-1)1,x)=\sum_{n\in \Z}a(n)x^{-n-1}
\]
for $a\in \fg$.  We have
\begin{eqnarray*}
\lefteqn{Y(a^{(0)}(-n_0)a^{(1)}(-n_1)\cdots a^{(r)}(-n_r)1,x)=}\\
&&=\res_{x_1}(x_1-x)^{-n_0}Y(a^{(0)}(-1)1,x_1)Y(a^{(1)}(-n_1)\cdots
a^{(r)}(-n_r)1,x)\\
&&\hspace{1em}-\res_{x_1}(-x+x_1)^{-n_0}Y(a^{(1)}(-n_1)\cdots
a^{(r)}(-n_r)1,x)Y(a^{(0)}(-1)1,x_1).\\
\end{eqnarray*}
Any restricted $\hat\fg$-module $W$ of level $\ell$ has a 
unique $V_{\hat\fg}(\ell,0)$-module
structure with the same action as above. Furthermore, in case
$\ell\neq -h$,
\[
\omega=\frac{1}{2(\ell+h)}\sum_{i=1}^{\dim\fg}g^i(-1)^2 1
\]
is a conformal element, where $\{g^i\}_{i=1,\dots,\dim\fg}$ is an
orthonormal basis of $\fg$ with respect to the form $(\cdot,\cdot)$,
and the quadruple $(V_{\hat\fg}(\ell,0), Y, {\bf 1}, \omega)$ is a
vertex operator algebra.
\end{theo}

Since every object of ${\cal O}_\kappa$ is a restricted
$\hat\fg$-module of level $\kappa-h$, it is a module for the vertex
algebra $V_{\hat\fg}(\kappa-h,0)$, and when $\kappa\neq 0$, it is a
generalized module for $V_{\hat\fg}(\kappa-h,0)$ viewed as a vertex
operator algebra.

For any element $\lambda\in W_1\circ_{Q(z)}W_2$ (recall
(\ref{circ_Q})), by Corollary \ref{coro} and the fact that $\fg\otimes
t^N\C[t,(z+t)^{-1}]\subset(\fg\otimes t\C[t,(z+t)^{-1}])^N$ in
$U(\fg\otimes\C((t)))$ we see that $\lambda$ satisfies the
$Q(z)$-compatibility condition. On the other hand, the closedness of
tensor product in ${\cal O}_\kappa$ in Theorem \ref{closethm} shows
that $(W_1\circ_{Q(z)}W_2)'$, and hence $W_1\circ_{Q(z)}W_2$ itself, is an
object of ${\cal O}_\kappa$. So all elements of $W_1\circ_{Q(z)}W_2$
also satisfy the $Q(z)$-local grading restriction condition. Hence
$W_1\circ_{Q(z)}W_2\subseteq W_1\hboxtr_{Q(z)}W_2$.

Conversely, if $\lambda\in W_1\hboxtr_{Q(z)}W_2$, then from the
$Q(z)$-local grading restriction condition we see that $\lambda$
generates a generalized $V$-module. But since for any $a\in\fg$,
$a\otimes t=(a(-1)1)_{1}$ as an operator on $W_1\hboxtr_{Q(z)}W_2$
reduces generalized weights by $1$, we see that when $N$ is large
enough we have $\xi_1\cdots\xi_N\lambda=0$ for all
$\xi_1,\dots,\xi_N\in\fg\otimes t\C[t,(z+t)^{-1}]$. Hence
$W_1\hboxtr_{Q(z)}W_2\subseteq W_1\circ_{Q(z)}W_2$.

We have proved:
\begin{theo}\label{qz-tr-prod}
For $W_1$ and $W_2$ in ${\cal O}_\kappa$, the two subspaces
$W_1\circ_{Q(z)}W_2$ and $W_1\hboxtr_{Q(z)}W_2$ of $(W_1\otimes
W_2)^*$ are equal to each other. In particular, the tensor product of
two modules constructed in \cite{KL2} with respect to $Q(z)$ satisfies
the universal property in Definition \ref{Qtpdef}.  \epf
\end{theo}

Now we proceed to the existence and construction of the associativity
isomorphism for this tensor product. For this, we now work in the
setting of tensor product associated with another type of Riemann
spheres with punctures and local coordinates, namely, the spheres with
punctures and local coordinates of type $P(z)$; recall from \cite{H1}
that for a nonzero complex number $z$, $P(z)$ denotes the Riemann
sphere with ordered punctures $\infty$, $z$, $0$ and local coordinates
$1/w$, $w-z$, $w$ around these punctures.

\begin{rema}\label{pz-tr-prod}
{\rm
The reason that we use $P(z)$ here is because it is most convenient
for the formulation of the associativity isomorphisms, due to the fact that
spheres with punctures and local coordinates of type $P(z)$ are closed
under sewing and subsequently decomposing. The corresponding
associativity isomorphisms for other types of tensor products can be
constructed from those for the type $P(z)$ by natural transformations
associated to certain parallel transport over the moduli space of
spheres with punctures and local coordinates. The $P(z)$-tensor product 
of  $W_1$ and $W_2$ in ${\cal O}_\kappa$ exists if and only if 
their $Q(z)$-tensor product exists; the details are given in \cite{HLZ2}.}
\end{rema}

We need the following notions from \cite{HLZ1} and \cite{HLZ2}, which
generalize the corresponding notions in \cite{tensor4} to the
logarithmic case:
\begin{defi}{\rm
Let $V$ be a vertex operator algebra and ${\cal C}$ be a category of
generalized $V$-modules. We say that products of intertwining
operators in ${\cal C}$ satisfy the {\em convergence and extension
property} if for any objects $W_1$, $W_2$, $W_3$, $W_4$ and $M_1$ of
${\cal C}$, and intertwining operator ${\cal Y}_1$ and ${\cal Y}_2$ of
types ${W_4\choose W_1\,M_1}$ and ${M_1\choose W_2\,W_3}$,
respectively, there exists an integer $N$ depending only on ${\cal
Y}_1$ and ${\cal Y}_2$, and for any $w_{(1)}\in W_1$, $w_{(2)}\in W_2$,
$w_{(3)}\in W_3$, $w'_{(4)}\in W'_4$, there exist $M\in{\mathbb N}$,
$r_{k}, s_{k}\in {\mathbb R}$, $i_{k}, j_{k}\in {\mathbb N}$, $k=1,\dots,M$
and analytic functions $f_{i_{k}j_{k}}(z)$ on $|z|<1$, $k=1, \dots,
M$, satisfying
\[
\wt w_{(1)}+\wt w_{(2)}+s_{k}>N, \;\;\;k=1, \dots, M,
\]
such that
\[
\langle w'_{(4)}, {\cal Y}_1(w_{(1)}, x_2) {\cal Y}_2(w_{(2)},
x_2)w_{(3)}\rangle_{W_4} \lbar_{x_1= z_1, \;x_2=z_2}
\]
is convergent when $|z_1|>|z_2|>0$ and can be analytically extended to
the multi-valued analytic function
\[
\sum_{k=1}^{M}z_2^{r_{k}}(z_1-z_2)^{s_{k}}(\log z_2)^{i_{k}}
(\log(z_1-z_2))^{j_{k}}f_{i_{k}j_{k}}\left(\frac{z_1-z_2}{z_2}\right)
\]
in the region $|z_2|>|z_1-z_2|>0$.
}
\end{defi}

\begin{defi}{\rm
Let $V$ be a vertex operator algebra and ${\cal C}$ be a category of
generalized $V$-modules. If for any $n\in{\mathbb Z}_+$, any generalized
$V$-modules $W_i$, $i=0, \dots, n+1$ and $\tilde W_i$, $i=1, \dots,
n-1$ in ${\cal C}$, and intertwining operators ${\cal Y}_{1}, {\cal
Y}_{2}, \dots, {\cal Y}_{n-1}, {\cal Y}_{n}$, of types ${W'_{0}\choose
W_{1}\widetilde{W}_{1}}$, ${\widetilde{W}_{1}\choose
W_{2}\widetilde{W}_{2}}, \dots, {\widetilde{W}_{n-2}\choose
W_{n-1}\widetilde{W}_{n-1}}$, ${\widetilde{W}_{n-1}\choose
W_{n}W_{n+1}}$, respectively, the series \begin{equation} \langle
w_{0}, {\cal Y}_{1}(w_{1}, z_{1})\cdots {\cal Y}_{n}(w_{n},
z_{n})w_{n+1}\rangle
\end{equation}
is absolutely convergent in the region $|z_{1}|>\cdots> |z_{n}|>0$,
then we say that {\em products of arbitrary number of intertwining
operators among objects of ${\cal C}$ are convergent}.  }
\end{defi}

The following theorem was proved in \cite{HLZ1} and \cite{HLZ2}:
\begin{theo}\label{conc}
Let $V$ be a vertex operator algebra of central charge $c\in\C$ and
${\cal C}$ a category of generalized $V$-modules closed under the
contragredient functor $(\cdot)'$ and under taking direct sums and
quotients.  Assume that the convergence and extension property for
products of intertwining operators holds in ${\cal C}$ and that products of
arbitrary number of intertwining operators among objects of ${\cal C}$
are convergent.  Further assume that every finitely-generated lower
truncated generalized $V$-module is an object of ${\cal C}$ and, for
any generalized $V$-modules $W_{1}$ and $W_{2}$,
$W_{1}\hboxtr_{P(z)}W_{2}$ is an object of ${\cal C}$.  Then the
category ${\cal C}$ has a natural structure of vertex tensor category
(see \cite{tensorK}) of central charge equal to the central charge $c$
of $V$ such that for each $z\in {\mathbb C}^{\times}$, the tensor product
bifunctor $\boxtimes_{\psi(P(z))}$ associated to $\psi(P(z))\in
\tilde{K}^{c}(2)$ is equal to $\boxtimes_{P(z)}$.  In particular, the
category ${\cal C}$ has a braided tensor category structure.
\end{theo}

By definition it is clear that ${\cal O}_\kappa$ is closed under
taking direct sums and quotients.  By Theorems \ref{closethm} and  
\ref{qz-tr-prod} and Remark \ref{pz-tr-prod}, we
have that ${\cal O}_\kappa$ is closed under the contragredient functor
and the $Q(z)$- and $P(z)$-tensor product functors. We now prove:
\begin{propo}
Any lower truncated, finitely-generated generalized
$V_{\hat\fg}(\ell,0)$-module is an object of ${\cal O}_\kappa$.
\end{propo}
\pf Let $W$ be a lower truncated generalized $V$-module generated by a
finite set $S$ of elements. That is, we have
\[
W=U(\hat\fg)S=U(\hat\fg_{(>)})U(\hat\fg_{(\leq)})S.
\]
Let $N=U(\fg_{(\leq)})S$, the $U(\fg_{(\leq)})$-submodule of $W$
generated by $S$. Then since $W$ is lower truncated and $S$ is finite,
$N$ is finite-dimensional. Let $N^\kappa$ be the $\hat\fg$-module
induced by $N$ as in (\ref{Nkappa}). Then $N^\kappa$ is a generalized
Weyl module and there is a unique $\hat\fg$-homomorphism from
$N^\kappa$ to $W$ fixing $N$. It is clear that this is a surjection
and thus $W$ is a quotient module of generalized Weyl module
$N^\kappa$. Hence $W$ is in ${\cal O}_\kappa$. \epf

\begin{rema}{\rm
Note that in general a generalized Weyl module may not be generated by
its lowest weight subspace, even if it is indecomposable. }
\end{rema}

Let $V$ be a vertex operator algebra. 
A generalized $V$-module $W$ is {\it $C_{1}$-cofinite} if 
$W/C_{1}(W)$ is finite dimensional, where $C_{1}(W)$ is the 
subspace of $W$ spanned by elements of the form 
$u_{-1}w$ for $u\in V_{+}=\coprod_{n\in \mathbb{N}}V_{(n)}$
and $w\in W$. 

\begin{propo}\label{c-1-cofinite}
The objects
of ${\cal O}_\kappa$ are $C_1$-cofinite as generalized 
$V_{\hat\fg}(\ell,0)$-modules. 
\end{propo}
\pf
Let $W$ be an object of ${\cal O}_\kappa$. Then 
$W$ is a quotient of a generalized Weyl module
$U(\hat\fg)\otimes_{U(\hat\fg_{(\leq 0)})}M$.
Let $\tilde{M}$ be the image of $M$ under projection 
{}from the generalized Weyl module to $W$. 
Then 
$W$ is spanned by elements of 
$\tilde{M}$ together with elements of the 
form $a(-n)w$ for $w\in W$ and $n\in \Z_{+}$. 
By the $L(-1)$-derivative property, we have 
\begin{eqnarray*}
a(-n)&=&(a(-1)\mathbf{1})_{-n}\nno\\
&=&\frac{(L(-1)a(-1)\mathbf{1})_{-n+1}}{n-1}
\end{eqnarray*}
when $n\ne 1$. 
Thus we see that $W$ is spanned by elements of 
$\tilde{M}$ together with elements of the form 
$u_{-1}w$ for $u\in (V_{\hat\fg}(\ell,0))_{+}$ and $w\in W$. 
Since $\tilde{M}$ is finite dimensional, we see that 
$W/C_{1}(W)$ is also finite dimensional. 
\epfv

Let $V$ be a vertex operator algebra. A generalized $V$-module 
$W$ is {\it quasi-finite dimensional} if 
for any real number $N$, $\coprod_{\Re{(n)}<N}W_{[n]}$ is
finite dimensional.

We have the following:

\begin{propo}\label{quasi-f-d}
The objects
of ${\cal O}_\kappa$ are quasi-finite dimensional as generalized 
$V_{\hat\fg}(\ell,0)$-modules. 
\end{propo}
\pf
Since generalized Weyl modules are obviously quasi-finite dimensional,
the objects of ${\cal O}_\kappa$, as quotients of 
generalized Weyl modules,
are quasi-finite dimensional.
\epfv

By Theorems \ref{closethm} and \ref{qz-tr-prod}, 
Remark \ref{pz-tr-prod}, Corollary \ref{c-1-cofinite}, Proposition 
\ref{quasi-f-d} and Theorem 7.2 in \cite{HLZ1}, we obtain the main conclusion 
of the present paper:

\begin{theo}
The category ${\cal O}_\kappa$ has a
natural structure of vertex tensor category and in particular,
a natural structure of braided tensor category.
\end{theo}

\bigskip

\noindent {\small \sc Department of Mathematics, Rutgers University,
Piscataway, NJ 08854}

\noindent {\em E-mail address}: linzhang@math.rutgers.edu

\end{document}